\documentclass[11pt,reqno]{amsart}
\usepackage{amssymb}
\usepackage[all]{xy}
\newtheorem{thm}{Theorem}[section]
\newtheorem{lem}[thm]{Lemma}
\newtheorem{prop}[thm]{Proposition}
\newtheorem{cor}[thm]{Corollary}

\theoremstyle{definition}
\newtheorem{dfn}[thm]{Definition}

\newtheorem{ex}[thm]{Example}

\theoremstyle{remark}
\newtheorem{remark}[thm]{Remark}

\newtheorem{notation}[thm]{Notation}

\newtheorem{assumption}[thm]{Assumption}
\newcommand{\CA}{{\mathcal{A}}}

\newcommand{\CE}{{\mathcal{E}}}
\newcommand{\bE}{{\overline{\mathcal{E}}}}

\newcommand{\CEa}{\CE^{0,-}}

\newcommand{\CL}{{\mathcal{L}}}
\newcommand{\CB}{{\mathcal{B}}}

\newcommand{\af}{\alpha}
\newcommand{\bt}{\beta}
\newcommand{\gm}{\gamma}
\newcommand{\dt}{\delta}

\newcommand{\ld}{\lambda}

\newcommand{\sm}{\sigma}

\newcommand{\Om}{\Omega}

\newcommand{\C}{{\mathbb{C}}}

\newcommand{\T}{{\mathbb{T}}}

\begin{document}

%\numberwithin{equation}{section}

\title[On simple labelled graph $C^*$-algebras]
{On simple labelled graph $C^*$-algebras}

\author[J. A. Jeong]{Ja A Jeong$^{\dagger}$$^{\ddagger}$}
\thanks{Research partially supported by KRF-2008-314-C00014$^{\dagger}$
and NRF-2009-0068619$^{\ddagger}$.}
\address{
Department of Mathematical Sciences and Research Institute of Mathematics\\
Seoul National University\\
Seoul, 151--747\\
Korea} \email{jajeong\-@\-snu.\-ac.\-kr }

\author[S. H. Kim]{Sun Ho Kim$^{\ddagger}$}
\address{
Department of Mathematical Sciences\\
Seoul National University\\
Seoul, 151--747\\
Korea} \email{hoya4200\-@\-snu.\-ac.\-kr}





\begin{abstract}
We consider the simplicity of the $C^*$-algebra associated to
a labelled space $(E,\CL,\bE)$, where
$(E,\CL)$ is a labelled graph and
$\bE$ is the smallest accommodating set
containing all generalized vertices.
We prove that if $C^*(E, \CL, \bE)$ is simple, then  $(E, \CL, \bE)$  is
strongly cofinal, and if, in addition,
$\{v\}\in \bE$ for every vertex $v$, then $(E, \CL, \bE)$ is disagreeable.
It is observed that $C^*(E, \CL, \bE)$ is simple
whenever  $(E, \CL, \bE)$ is strongly cofinal and disagreeable,
which is recently known for the $C^*$-algebra $C^*(E, \CL, \CEa)$
associated to a labelled space $(E, \CL, \CEa)$ of
the smallest accommodating set $\CEa$.
\end{abstract}

\maketitle

\setcounter{equation}{0}

%\section{Introduction}

\section{Introduction}
Given a directed graph $E=(E^0,E^1)$ with the vertex set
$E^0$ and the edge set $E^1$, the  $C^*$-algebra $C^*(E)$
generated by the family of universal Cuntz-Krieger $E$-family is
associated (see \cite{KPR, KPRR, Ra} among others).
A Cuntz-Kriger algebra can
now be viewed as a graph $C^*$-algebra  of a finite graph.
It is well known that a graph $C^*$-algebra  of a row finite graph $E$
is simple if and only if the graph $E$ is cofinal
and every loop has an exit (\cite{KPR}).

Let  $E$ be a directed graph and $\CA$ be a countable alphabet.
If  $\CL:E^1\to \CA$ is a labelling map, we call $(E,\CL)$
a labelled graph.
If $\CB\subset 2^{E^0}$ is an accommodating set for $(E,\CL)$,
there exists a labelled graph  $C^*$-algebra $C^*(E,\CL,\CB)$
associated to the labelled space $(E,\CL,\CB)$ (\cite{BP1}).
 Graph $C^*$-algebras and more generally 
 ultragraph $C^*$-algebras (\cite{To1,To2})
 are labelled graph $C^*$-algebras (\cite{BP1}).

The simplicity of a labelled graph $C^*$-algebra
$C^*(E,\CL,\CB)$ is considered in
\cite{BP2} when  $\CB$ is the smallest accommodating set $\CEa$.
Using the generalized vertices that play the role of the vertices in
a directed graph one has the useful expression of
the elements of a dense set in $C^*(E,\CL,\CEa)$.
But a generalized vertex does not necessarily belong to
the accommodating set $\CEa$.
In this paper we consider the simplicity of the labelled graph $C^*$-algebra
$C^*(E,\CL,\bE)$ where $\bE$ is the smallest accommodating set
containing the generalized vertices.

We first review in the next section the definition of a graph $C^*$-algebra
(from \cite{KPR, KPRR}) and
a labelled graph $C^*$-algebra (from \cite{BP1, BP2}) to set up the notations.
Then we show in section 3 that if a labelled graph $C^*$-algebra
$C^*(E,\CL,\bE)$ is simple, the labelled space $(E,\CL,\bE)$ is
strongly cofinal (Theorem~\ref{thm-stcofinal}).
If, in addition,  $\{v\}\in \bE$ for each $v\in E^0$, it is shown that
$(E,\CL,\bE)$ is disagreeable (Theorem~\ref{thm-disagreeable}).
The fact that if a labelled space
$(E,\CL,\bE)$ is strongly cofinal and disagreeable then
$C^*(E,\CL,\bE)$ is simple can be obtained by a slight modification
of the proof of \cite[Theorem 6.4]{BP2} (Theorem~\ref{thm-simple}).

\vskip 1pc

\section{Labelled spaces and their $C^*$-algebras $C^*(E, \CL, \CB)$ }

 A  {\it directed graph} $E=(E^0,E^1,r,s)$ consists of the vertex set $E^0$,
 the edge set $E^1$, and the range, source maps $r_E,s_E: E^1\to E^0$.
 We shall simply refer to  directed graphs as  graphs
 and often write $r$, $s$ for $r_E$, $s_E$.
 By $E^n$ we denote the set of all finite paths $\ld=\ld_1\cdots \ld_n$
 of {\it length} $n$ ($|\ld|=n$),
($\ld_{i}\in E^1,\ r(\ld_{i})=s(\ld_{i+1}), 1\leq i\leq n-1$).
We also set $E^{\leq n}:=\cup_{i=1}^n E^i$ and
$E^{\geq n}:=\cup_{i=n}^\infty E^i$.
The maps $r$ and $s$ naturally extend to $E^{\geq 1}$.
Infinite paths
$\ld_1\ld_2\ld_3\cdots $  can also
be considered $(r(\ld_{i})=s(\ld_{i+1}),\, \ld_i\in E^1)$
and  we denote the set of all infinite paths by
$E^\infty$ extending $s$  to $E^\infty$ by
$s(\ld_1\ld_2\ld_3\cdots):=s(\ld_1)$.

 A {\it labelled graph} $(E,\CL)$ over a countable alphabet set $\CA$
 consists of a directed graph $E$ and
 a {\it labelling map} $\CL:E^1\to \CA$.
 We assume that $\CL$ is onto.
 For a finite path $\ld=\ld_1\cdots \ld_n\in E^n$,
 put $\CL(\ld)=\CL(\ld_1)\cdots \CL(\ld_n)$.
 Similarly, we put
 $\CL(E^\infty)=\{\CL(\ld_1)\CL(\ld_2)\cdots\mid \ld_1\ld_2\cdots\in E^\infty\}$.
For a  path $\af=\af_1\cdots \af_n\in E^{\geq 1}\cup \CL(E^{\geq 1})$
and $1\leq k\leq l\leq n<\infty$,
we set $\af_{[k,l]}:=\af_k\cdots \af_l$
(and $\af_{[k,\infty)}:=\af_k\af_{k+1}\cdots$
for $\af\in E^\infty\cup \CL(E^\infty)$).
 The {\it range} and {\it source} of a labelled path $\af\in \CL(E^{\geq 1})$
are subsets of $E^0$ defined by
 $$r_{\CL_E}(\af)=\{r(\ld) \mid \ld\in E^{\geq 1},\,\CL(\ld)=\af\},\
 s_{\CL_E}(\af)=\{s(\ld) \mid \ld\in E^{\geq 1},\, \CL(\ld)=\af\}.$$
 (For both a labelled graph  and its underlying graph,
 we usually use the same notation $r$ and $s$ to denote range maps
 $r_E$, $r_{\CL_E}$ and source maps $s_E$, $s_{\CL_E}$.)
 The {\it relative range of $\af\in \CL(E^{\geq 1})$
 with respect to $A\subset 2^{E^0}$} is defined to be
$$
 r(A,\af)=\{r(\ld)\mid \ld\in E^{\geq 1},\ \CL(\ld)=\af,\ s(\ld)\in A\}.
$$
 If $\CB\subset 2^{E^0}$ is a collection of subsets of $E^0$ such that
 $r(A,\af)\in \CB$ whenever $A\in \CB$ and $\af\in \CL(E^{\geq 1})$,
 $\CB$ is said to be
 {\it closed under relative ranges} for $(E,\CL)$.
 We call $\CB$ an {\it accommodating set} for $(E,\CL)$
 if it is closed under relative ranges,
 finite intersections and unions and
 contains $r(\af)$ for all $\af\in \CL(E^{\geq 1})$.
  If $\CB$ is accommodating for $(E,\CL)$, the triple $(E,\CL,\CB)$ is called
 a {\it labelled space}.
 A labelled space $(E,\CL,\CB)$ is {\it weakly left-resolving} if
 $$r(A,\af)\cap r(B,\af)=r(A\cap B,\af)$$
 for every $A,B\in \CB$ and every $\af\in \CL(E^{\geq 1})$.
 If $(E,\CL,\CB)$ is weakly left-resolving
 and  $A,\, B,\, A\setminus B\in \CB$, then
 it follows that for $\af\in \CL(E^{\geq 1})$
 \begin{eqnarray}\label{differenceset}
 r(A\setminus B, \af)=r(A,\af)\setminus r(B,\af).
 \end{eqnarray}
 But (\ref{differenceset}) may not hold if $A\setminus B\notin \CB$
 as we see in Example~\ref{ex-rangevertex}.
 For $A,B\in 2^{E^0}$ and $n\geq 1$, let
 $$ AE^n =\{\ld\in E^n\mid s(\ld)\in A\},\ \
  E^nB=\{\ld\in E^n\mid r(\ld)\in B\},$$
 and
 $AE^nB=AE^n\cap E^nB$.
 We write $E^n v$ for $E^n\{v\}$ and $vE^n$ for $\{v\}E^n$,
 and will use notations like $AE^{\geq k}$ and $vE^\infty$
 which should have obvious meaning.
 A labelled graph $(E,\CL)$ is {\it left-resolving}
if for all $v\in E^0$ the map
$\CL:E^1 v\to \CL(E^1v)$ is bijective
(hence $\CL: E^{\leq k}v \to \CL(E^{\leq k}v)$ is bijective
for all $k\geq 1$).
 A labelled space $(E,\CL,\CB)$ is said to be {\it set-finite}
 ({\it receiver set-finite}, respectively) if for every $A\in \CB$
 the set  $\CL(AE^1)$ ($\CL(E^1A)$, respectively) finite.

\vskip 1pc

\begin{dfn} (\cite[Definition 4.1]{BP1})
\label{definition-representation}
Let $(E,\CL,\CB)$ be a weakly left-resolving labelled space.
A {\it representation} of $(E,\CL,\CB)$
consists of projections $\{p_A: A\in \CB\}$ and
partial isometries
$\{s_a: a\in \CA\}$ such that for $A, B\in \CB$ and $a, b\in \CA$,
\begin{enumerate}
\item[(i)] $p_{A\cap B}=p_Ap_B$ and $p_{A\cup B}=p_A+p_B-p_{A\cap B}$,
where $p_{\emptyset}=0$,
\item[(ii)] $p_A s_a=s_a p_{r(A,a)}$,
\item[(iii)] $s_a^*s_a=p_{r(a)}$ and $s_a^* s_b=0$ unless $a=b$,
\item[(iv)] for $A\in \CB$, if  $\CL(AE^1)$ is finite and non-empty, then
$$p_A=\sum_{a\in \CL(AE^1)} s_a p_{r(A,a)}s_a^*.$$
\end{enumerate}

\noindent It is known \cite[Theorem 4.5]{BP1} that
if $(E,\CL,\CB)$ is a weakly left-resolving labelled space,
there exists a $C^*$-algebra $C^*(E,\CL,\CB)$
generated by a universal representation
$\{s_a,p_A\}$ of $(E,\CL,\CB)$.
We call $C^*(E,\CL,\CB)$ the {\it labelled graph $C^*$-algebra} of
a labelled space $(E,\CL,\CB)$.
\end{dfn}

\vskip 1pc

\begin{assumption}\label{assumptions-graph}
Throughout this paper, we assume the following.
\begin{enumerate}
\item[(a)]  $E$ has no sinks, that is
$|s^{-1}(v)|>0$ for all $v\in E^0$, and
$(E,\CL,\CB)$ is set-finite and receiver set-finite.

\item[(b)]
If $e,f\in E^1$ are edges with $s(e)=s(f)$, $r(e)=r(f)$, and
$\CL(e)=\CL(f)$, then $e=f$.
\end{enumerate}
\end{assumption}

\vskip 1pc

\begin{remark}\label{elements}
Let $C^*(E,\CL,\CB)$ be the labelled graph $C^*$-algebra of
$(E,\CL,\CB)$ with generators $\{s_a,p_A\}$.
Then  $s_a\neq 0$ and $p_A\neq 0$ for $a\in \CA$
and  $A\in \CB$, $A\neq \emptyset$.
Note also that
$s_\af p_A s_\bt^*\neq 0$ if and only if $A\cap r(\af)\cap r(\bt)\neq\emptyset$.
Since we assume that $(E,\CL,\CB)$ is set-finite and $E$ has no sinks,
by \cite[Lemma 4.4]{BP1} and
Definition~\ref{definition-representation}(iv)
it follows that
$$p_A=\sum_{\sm\in \CL(AE^n)} s_\sm p_{r(A,\sm)}s_\sm^*
\ \text{ for }   A\in \CB \text{ and } n\geq 1$$
and
$$C^*(E,\CL,\CB)=\overline{span}\{s_\af p_A s_\bt^*\mid
\af,\,\bt\in \CL(E^{\geq 1}),\ A\in \CB\}.$$
\end{remark}

\vskip 1pc

\noindent

For $v,w\in E^0$, we write $v\sim_l w$ if $\CL(E^{\leq l} v)=\CL(E^{\leq l} w)$
as in \cite{BP2}.
Then  $ \sim_l $ is an equivalence relation on $E^0$.
The equivalence class $[v]_l$ of $v$  is called a {\it generalized vertex}.
Let $\Om_l:=E^0/\thicksim_l$.
   For $k<l$ and $v\in E^0$,  $[v]_l\subset [v]_k$ is obvious and
   $[v]_l=\cup_{i=1}^m [v_i]_{l+1}$
   for some vertices  $v_1, \dots, v_m\in [v]_l$ (\cite[Proposition 2.4]{BP2}).

Let $\CEa$ be the smallest accommodating set for $(E,\CL)$.
Then
$$\CEa=\{ \cup_{k=1}^m\cap_{i=1}^n  r(\bt_{i,k})\mid \bt_{i,k}\in \CL(E^{\geq 1})\}$$
(see \cite[Remark 2.1]{BP2}).
Also every accommodating set $\CB$ for $(E,\CL)$ contains $\CEa$.
By  \cite[Proposition 2.4]{BP2},
every generalized vertex $[v]_l$ is the difference of two sets
$X_l(v)$ and $r(Y_l(v))$ in $\CEa$,
more precisely,
$$[v]_l=X_l(v)\setminus r(Y_l(v)),$$  where
$X_l(v):=\cap_{\af\in \CL(E^{\leq l} v)} r(\af)\ \text{ and }\
 Y_l(v):=\cup_{w\in X_l(v)} \CL(E^{\leq l} w)\setminus \CL(E^{\leq l} v).$
The following example shows that
$$r([v]_l, a) \neq r(X_l(v),a)\setminus r(r(Y_l(v)), a),$$
in general.

\vskip 1pc

\begin{ex}\label{ex-rangevertex}
Consider the following weakly left-resolving labelled space
$(E, \CL,\CEa)$:

\vskip 1pc

\hskip .5pc
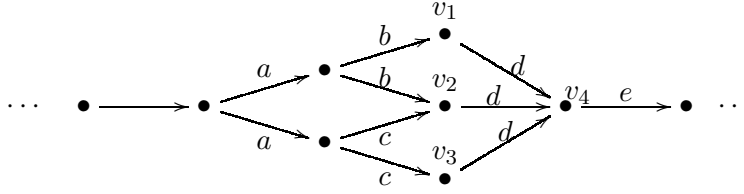
\xy
/r0.38pc/:(-15,0)*+{\cdots}="V-1";
(-10,0)*+{\bullet}="V0";
(0,0)*+{\bullet}="V1";
 "V0";"V1"**\crv{(-10,0)&(-5,0)&(0,0)};
 ?>*\dir{>}\POS?(.5)*+!D{};
 (10,3)*+{\bullet}="V2u";
 (10,-3)*+{\bullet}="V2l";
 "V1";"V2u"**\crv{(0,0)&(10,3)};?>*\dir{>}\POS?(.5)*+!D{ };
 "V1";"V2l"**\crv{(0,0)&(10,-3)};?>*\dir{>}\POS?(.5)*+!D{ };
 (20,0)*+{\bullet}="V30";
 (20,6)*+{\bullet}="V3u";
 (20,-6)*+{\bullet}="V3l";
 "V2u";"V30"**\crv{(10,3)&(20,0)};?>*\dir{>}\POS?(.5)*+!D{},
 "V2l";"V30"**\crv{(10,-3)&(20,0)};?>*\dir{>}\POS?(.5)*+!D{},
 "V2u";"V3u"**\crv{(10,3)&(20,6)};?>*\dir{>}\POS?(.5)*+!D{},
 "V2l";"V3l"**\crv{(10,-3)&(20,-6)};?>*\dir{>}\POS?(.5)*+!D{},
 (30,0)*+{\bullet}="V4";
 "V30";"V4"**\crv{(20,0)&(30,0)};?>*\dir{>}\POS?(.5)*+!D{ };
 "V3u";"V4"**\crv{(20,6)&(30,0)};?>*\dir{>}\POS?(.5)*+!D{ };
 "V3l";"V4"**\crv{(20,-6)&(30,0)};?>*\dir{>}\POS?(.5)*+!D{ };
 (40,0)*+{\bullet}="V5";
 "V4";"V5"**\crv{(30,0)&(40,0)};?>*\dir{>}\POS?(.5)*+!D{ };
 (5,3)*+{a},(5,-3)*+{a},(15,6)*+{b},(15,2.5)*+{b},
 (15,-2.6)*+{c},(15,-6)*+{c},(26,3.5)*+{d},(24,1)*+{d},(25,-2)*+{d},(35,1)*+{e},
 (20,8)*+{v_1},(20,2)*+{v_2},(20,-4)*+{v_3},(31,1)*+{v_4}, (44,0)*+{\cdots};
(40,0)*+{ }="V6";
\endxy

\vskip 2pc

\noindent Since
$\{v_1,v_2\}=r(ab)$, $\{v_2, v_3\}=r(ac)$, and $\{v_2\}=\{v_1,v_2\}\cap \{v_2, v_3\}$,
we have $\{v_1,v_2\}, \,\{v_2, v_3\},\, \{v_2\}\in \CE^{0,-}$.
 From $X_2(v_1)=\{v_1,v_2\}, \ Y_2(v_1)=\{c,ac\}$ and $r(Y_2(v_1))=\{v_2,v_3\}$,
 we have
$$r(X_2(v_1),d)\setminus r(r(Y_2(v_1)),d )=\{v_4\}\setminus \{v_4\}=\emptyset.$$
But  $X_2(v_1)\setminus r(Y_2(v_1))=\{v_1\}=[v_1]_2 \notin \CE^{0,-}$
and $r([v_1]_2,d)=r(\{v_1\},d)=\{v_4\}\neq \emptyset$.

\end{ex}

\vskip 1pc

\section{Simplicity of a labelled graph  $C^*$-algebra  $C^*(E, \CL, \bE)$ }

\vskip 1pc

\noindent
   Recall (\cite{BP2}) that
   a labelled space $(E, \CL, \CB)$ is $l$-{\it cofinal} if for all
   $x\in  \CL(E^\infty)$, $[v]_l\in \Om_l$, and $w\in s(x)$,
   there are $R(w)\geq l$, $N\geq 1$, and  finitely many labelled paths
   $\ld_1, \dots, \ld_m$ such that for all $d\geq R(w)$ we have
$$r([w]_d,  x_{[1,N]})\subseteq \cup_{i=1}^m r([v]_l,\ld_i)$$
  Since $[w]_{d'}\subset [w]_{d}$ whenever $d\leq d'$,
  we may restate the  definition of an $l$-cofinal labelled space
  as follows.

\vskip 1pc
\begin{dfn} Let $(E, \CL, \CB)$ be a labelled space.
\begin{enumerate}
\item[(a)] (\cite{BP2}) $(E, \CL, \CB)$ is $l$-cofinal
 if for all  $x\in  \CL(E^\infty)$, $[v]_l\in \Om_l$, and $w\in s(x)$,
 there are $d\geq l$,  $N\geq 1$, and a finite number of labelled paths
 $\ld_1, \dots, \ld_m$ such that
$$r([w]_d, x_{[1,N]})\subseteq \cup_{i=1}^m r([v]_l,\ld_i).$$
$(E, \CL, \CB)$ is  {\it cofinal} if there is an $L>0$ such that
$(E, \CL, \CB)$ is $l$-cofinal  for all $l\geq L$.

\item[(b)] $(E, \CL, \CB)$ is  {\it strongly cofinal} if
for all  $x\in  \CL(E^\infty)$, $[v]_l\in \Om_l$,
and $w\in s(x)$,
there are $N\geq 1$ and  a finite number of labelled paths
$\ld_1, \dots, \ld_m$  such that
$$r([w]_1,  x_{[1,N]})\subset \cup_{i=1}^m r([v]_l, \ld_i).$$
\end{enumerate}

\end{dfn}

\vskip 1pc
\noindent
Every  strongly cofinal labelled space  is cofinal.
 Note that $(E, \CL, \CB )$ is  not  cofinal
 if and only if there is a sequence $l_1<l_2<\cdots$
 of positive integers such that $(E, \CL, \CB)$ is  not $l_i$-cofinal for all   $i\geq 1$.
 But, in fact, we have the following.

\vskip 1pc

\begin{prop}  A labelled space $(E, \CL, \CB)$ is  cofinal if and only if it is  $l$-cofinal
 for all $l\geq1$.
\end{prop}

\begin{proof} If the labelled space is $l$-cofinal for all $l\geq1$, obviously it is cofinal.
For the converse, it suffices to show that if $(E, \CL, \CB)$ is not $l$-cofinal,
it is not $l'$-cofinal for all $l'>l$.
To see this, first note that $(E, \CL, \CB)$ is  not $l$-cofinal if and only if
the set $\triangle_l$ of  triples $(x,[v]_l, w)$ of
$x\in  \CL(E^\infty)$, $[v]_l\in \Om_l$, and $w\in s(x)$ such that
$$r([w]_d,  x_{[1,N]})\nsubseteq \cup_{i=1}^m r([v]_l,\ld_i) $$
for all $d\geq l$, $N\geq 1$, and
a finite number of labelled paths $\ld_1, \dots, \ld_m$ is nonempty.
 Suppose $(E, \CL, \CB)$ is not $l$-cofinal and $(x,[v]_l, w)\in \triangle_{l}$.
 If $(E, \CL, \CB)$ is $l'$-cofinal for $l'>l$,
 then there are $d\geq l'$, $N\geq 1$, and
 a finite number of labelled paths $\ld_1, \dots, \ld_m$ such that
 $$r([w]_d,  x_{[1,N]})\subseteq \cup_{i}r([v]_{l'},\ld_i).$$
But this  contradicts to $(x,[v]_l, w)\in \triangle_l$
since $\cup_{i}r([v]_{l'},\ld_i)\subseteq \cup_{i}r([v]_{l},\ld_i)$.
\end{proof}

\vskip 1pc

\begin{notation}
 Let $\bE$  denote the set of all finite unions of generalized vertices.
 Then  by \cite[Proposition 2.4]{BP2},
 $$\CEa\subset \bE=\{\cup_{i=1}^m [v_i]_{l}\mid v_i\in E^0,\ m,\, l\geq 1\}.$$

\end{notation}

\vskip 1pc

 For the proof of the following proposition it is helpful to note that
 if $[v]_l\cap [w]_k\neq\emptyset$,
 then either $[v]_l\subset [w]_k$ or $[w]_k\subset [v]_l$.

\vskip 1pc

\begin{prop}\label{prop-cE}
Let $(E,\CL)$ be a labelled graph such that
\begin{eqnarray}\label{weaklr}  r(A,\af)\cap r(B,\af)=r(A\cap B,\af)
\end{eqnarray}
for all $A,B\in \bE$ and   $\af\in \CL(E^{\geq 1})$.
Then $\bE$ is a weakly left-resolving accommodating set such that
$A\setminus B\in \bE$  for $A, B\in \bE$.
If  $(E,\CL,\CEa)$ is weakly left-resolving and
$[v]_l\in \CEa$ for all $v\in E^0$ and $l\geq 1$,
then $\bE=\CEa$.
\end{prop}
\begin{proof}
Clearly $\bE$ is closed under finite intersections and unions.
    We show that $\bE$ is closed under relative ranges.
   For this, note first that
   $A\setminus B\in \bE,\ \ A, B\in  \CEa$
 (\cite[Proposition 2.4]{BP2}).
  Since the two sets  $[v]_l$ and $r(Y_l(v)) (\in \CEa)$ are disjoint,
  if $s(\af)\cap [v]_l\neq \emptyset$ and $s(\af)\cap  r(Y_l(v))\neq \emptyset$,
  then  by   (\ref{weaklr}),
    $$r([v]_l,\af)\cap r(r(Y_l(v)),\af)=r([v]_l\cap r(Y_l(v)),\af)=\emptyset.$$
  Hence $r([v]_l,\af) = r(X_l(v),\af)\setminus r(r(Y_l(v)),\af)\in \bE$
  since both  $r(X_l(v),\af)$ and $ r(r(Y_l(v)),\af)$ belong to $\CEa$.
  For an arbitrary $C=\cup_{i=1}^n [v_i]_l\in \bE$,
$$r(C,\af) =r(\cup_{i=1}^n [v_i]_l, \af)=\cup_{i=1}^n r([v_i]_l,\af)\in \bE.$$
\end{proof}

\vskip 1pc

\noindent
The labelled space  $(E, \CL, \CEa)$ in Example~\ref{ex-rangevertex}
is  weakly left-resolving, but $(E, \CL, \bE)$
is not weakly left-resolving
since for $\{v_1\}=[v_1]_2,\in \bE$,  $\{v_3\}=[v_3]_2 \in \bE$  we have
  $$r(\{v_1\} \cap \{v_3\}, d)=\emptyset \neq \{v_4\}=r(\{v_1\},d)\cap r(\{v_3\},d).$$
  We will consider only weakly left-resolving labelled spaces $(E, \CL, \bE)$
  for the rest of this paper,
  so Example~\ref{ex-rangevertex} is excluded from our discussion.

\vskip 1pc

Recall that a graph $E$ is {\it locally finite} if every vertex
receives and emits only finite number of edges.

\vskip 1pc

\begin{prop}\label{prop-nonsimple}
Let $E$ be a locally finite graph and
 $(E, \CL, \bE)$ be a weakly left-resolving labelled space
such that  $|r(a)|=\infty$ for an $a\in \CA$.
Suppose that for each $v\in E^0$, $[v]_l$ is finite for some $l\geq 1$.
Then  the projection $p_{[v]_l}$ of a finite set $[v]_l$
 generates a proper ideal of $C^*(E, \CL, \bE)$.
\end{prop}

\begin{proof}
Let  $I$ be the ideal
of $C^*(E, \CL, \bE)$ generated by the projection  $p_{[v]_l}$.
Suppose $I=C^*(E, \CL, \bE)$.
Then there exists an $X \in I$ such that $\|s_a^* s_a  - X\|<1$.
Let $X=\sum_{i=1}^m c_i(s_{\af_i} p_{A_i} s_{\bt_i}^*)p_{[v]_l}
(s_{\gm_i} p_{B_i} s_{\dt_i}^*)$,  $c_i\in\C$.
Then
$$V=\cup_{i=1}^m r([v]_l,\bt_i)\,\cup\, \cup_{i=1}^m r([v]_l, \gm_i)$$
is a finite set since  $E$ is locally finite, and $V \in \bE$.
Hence  the set
$$W:=\{s(\sm_i)\mid \CL(\sm_i)=\af_i,\ r(\sm_i)\cap V\neq\emptyset,\
\sm_i\in E^{\geq 1},\ i=1,2,\dots,m\}$$
is also finite.
Since  $r(a)\setminus V\,(\in \bE)$ is an infinite set,
one can choose $w\in r(a)\setminus V$
and $k\geq 1$ such that $[w]_k\cap W=\emptyset$.
Then $r([w]_k,\af_i) \cap r([v]_l,\bt_j)=\emptyset$
for all $i,j=1,\dots, m$ and
$p_{[w]_k}\leq p_{r(a)}=s_a^* s_a$.
Thus
$$ p_{[w]_k}(s_{\af_i} p_{A_i} s_{\bt_i}^*)p_{[v]_l}(s_{\gm_i} p_{B_i} s_{\dt_i}^*)
=  s_{\af_i} p_{r([w]_k,\af_i)\cap A_i\cap r([v]_l,\bt_i)}s_{\bt_i}^*(s_{\gm_i} p_{B_i} s_{\dt_i}^*)
= 0 $$
for all $i$, and then
$$
1>  \ \|s_a^* s_a - X\|\geq \|p_{[w]_k}(s_a^* s_a - X)p_{[w]_k}\|=\|p_{[w]_k}\|=1,$$
a contradiction.
\end{proof}

\vskip 1pc

\begin{ex}\label{ex-nonsimple}
Consider the following  labelled graph $(E, \CL)$ of \cite[7.2]{BP2}:

\vskip 1pc

\hskip .5pc
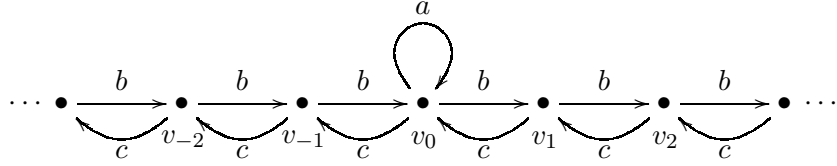
\xy
/r0.38pc/:(-33,0)*+{\cdots};
(33,0)*+{\cdots};
(-30,0)*+{\bullet}="V-3";
(-20,0)*+{\bullet}="V-2";
(-10,0)*+{\bullet}="V-1";
(0,0)*+{\bullet}="V0";
(10,0)*+{\bullet}="V1";
(20,0)*+{\bullet}="V2";
(30,0)*+{\bullet}="V3";
 "V-3";"V-2"**\crv{(-30,0)&(-20,0)};
 ?>*\dir{>}\POS?(.5)*+!D{};
 "V-2";"V-1"**\crv{(-20,0)&(-10,0)};
 ?>*\dir{>}\POS?(.5)*+!D{};
 "V-1";"V0"**\crv{(-10,0)&(0,0)};
 ?>*\dir{>}\POS?(.5)*+!D{};
 "V0";"V1"**\crv{(0,0)&(10,0)};
 ?>*\dir{>}\POS?(.5)*+!D{};
 "V1";"V2"**\crv{(10,0)&(20,0)};
 ?>*\dir{>}\POS?(.5)*+!D{};
 "V2";"V3"**\crv{(20,0)&(30,0)};
 ?>*\dir{>}\POS?(.5)*+!D{};
 "V-2";"V-3"**\crv{(-20,0)&(-25,-6)&(-30,0)};
 ?>*\dir{>}\POS?(.5)*+!D{};
 "V-1";"V-2"**\crv{(-10,0)&(-15,-6)&(-20,0)};
 ?>*\dir{>}\POS?(.5)*+!D{};
 "V0";"V-1"**\crv{(0,0)&(-5,-6)&(-10,0)};
 ?>*\dir{>}\POS?(.5)*+!D{};
 "V1";"V0"**\crv{(10,0)&(5,-6)&(0,0)};
 ?>*\dir{>}\POS?(.5)*+!D{};
 "V2";"V1"**\crv{(20,0)&(15,-6)&(10,0)};
 ?>*\dir{>}\POS?(.5)*+!D{};
 "V3";"V2"**\crv{(30,0)&(25,-6)&(20,0)};
 ?>*\dir{>}\POS?(.5)*+!D{};
 "V0";"V0"**\crv{(0,0)&(-4,4)&(0,8)&(4,4)&(0,0)};
 ?>*\dir{>}\POS?(.5)*+!D{};
 (-25,2)*+{b};(-15,2)*+{b};(-5,2)*+{b};(5,2)*+{b};(15,2)*+{b};(25,2)*+{b};
 (-25,-4)*+{c};(-15,-4)*+{c};(-5,-4)*+{c};(5,-4)*+{c};(15,-4)*+{c};(25,-4)*+{c};
 (0,8)*+{a};(0.1,-3)*+{v_0};(10.1,-3)*+{v_1};
 (-9.9,-3)*+{v_{-1}};
 (-19.9,-3)*+{v_{-2}};
 (20.1,-3)*+{v_{2}};
\endxy

\vskip 1pc
\noindent For each vertex $v_k$, we have
$ \{v_k\}=r(ab^k)$ for $k\geq 0$ and $ \{v_k\}=r(ac^k)$ for $k< 0$.
Hence $[v_k]_{k+1}=\{v_k\}\in \CEa$ for all $k$.
Despite the fact that $\CEa=\{A\subset E^0: A=E^0 \text{ or } A\text{ is a finite set}\}$
and  $\bE=\CEa\cup \{B\subset E^0: E^0\setminus B \text{ is a finite set}\}$,
one can show that $C^*(E,\CL,\CEa)\cong C^*(E,\CL,\bE)$
applying universal property of $C^*(E,\CL,\CEa)$ (since
$C^*(E,\CL,\bE)$ contains a family of generators that is a
representation of $C^*(E,\CL,\CEa)$)
and the gauge invariant uniqueness theorem
(\cite[Theorem 5.3]{BP1}).
 Since $|r(b)|=\infty$,
 by Proposition~\ref{prop-nonsimple},
 $C^*(E,\CL, \CEa)\cong C^*(E,\CL, \bE)$ is not simple
 (see Remark~\ref{remark-BP2}).
 But $(E,\CL,\CEa)$ is cofinal and
 disagreeable (we will discuss disagreeable labelled spaces later).
 Note that $(E,\CL,\CEa)$ is not strongly cofinal.
 In fact, for $[v_1]_1=\{v_k\mid k\neq 0\}$, $[v_{0}]_1=\{v_{0}\}$,
 and $x=bbb\cdots\in \CL(E^\infty)$, the set
 $r([v_1]_1,  x_{[1,N]} )$ is infinite  for any $N\geq  1$
 while $\cup_{i=1}^m r([v_0]_1, \ld_i)$ is finite
 for any finite number of  labelled paths  $\ld_1, \dots, \ld_m$.
 Hence it is not possible to have
 $r([v_1]_1,  x_{[1,N]} )\subset\cup_{i=1}^m r([v_0]_1, \ld_i)$.
\end{ex}

\vskip 1pc

\begin{remark}\label{remark-nonsimple}
The ideal $I$ generated by the projection $p_{\{v_0\}}$ in
Example~\ref{ex-nonsimple} is isomorphic to the graph $C^*$-algebra
$C^*(E)$, where $E$ is the following graph.
\vskip 1pc
\hskip .5pc
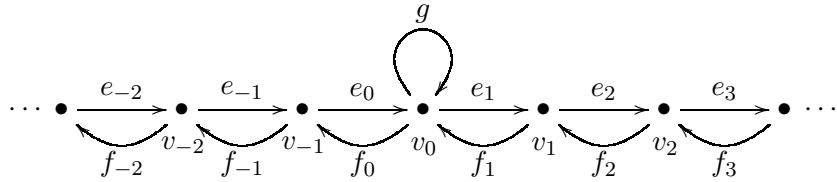
\xy
/r0.38pc/:(-33,0)*+{\cdots};
(33,0)*+{\cdots};
(-30,0)*+{\bullet}="V-3";
(-20,0)*+{\bullet}="V-2";
(-10,0)*+{\bullet}="V-1";
(0,0)*+{\bullet}="V0";
(10,0)*+{\bullet}="V1";
(20,0)*+{\bullet}="V2";
(30,0)*+{\bullet}="V3";
 "V-3";"V-2"**\crv{(-30,0)&(-20,0)};
 ?>*\dir{>}\POS?(.5)*+!D{};
 "V-2";"V-1"**\crv{(-20,0)&(-10,0)};
 ?>*\dir{>}\POS?(.5)*+!D{};
 "V-1";"V0"**\crv{(-10,0)&(0,0)};
 ?>*\dir{>}\POS?(.5)*+!D{};
 "V0";"V1"**\crv{(0,0)&(10,0)};
 ?>*\dir{>}\POS?(.5)*+!D{};
 "V1";"V2"**\crv{(10,0)&(20,0)};
 ?>*\dir{>}\POS?(.5)*+!D{};
 "V2";"V3"**\crv{(20,0)&(30,0)};
 ?>*\dir{>}\POS?(.5)*+!D{};
 "V-2";"V-3"**\crv{(-20,0)&(-25,-6)&(-30,0)};
 ?>*\dir{>}\POS?(.5)*+!D{};
 "V-1";"V-2"**\crv{(-10,0)&(-15,-6)&(-20,0)};
 ?>*\dir{>}\POS?(.5)*+!D{};
 "V0";"V-1"**\crv{(0,0)&(-5,-6)&(-10,0)};
 ?>*\dir{>}\POS?(.5)*+!D{};
 "V1";"V0"**\crv{(10,0)&(5,-6)&(0,0)};
 ?>*\dir{>}\POS?(.5)*+!D{};
 "V2";"V1"**\crv{(20,0)&(15,-6)&(10,0)};
 ?>*\dir{>}\POS?(.5)*+!D{};
 "V3";"V2"**\crv{(30,0)&(25,-6)&(20,0)};
 ?>*\dir{>}\POS?(.5)*+!D{};
 "V0";"V0"**\crv{(0,0)&(-4,4)&(0,8)&(4,4)&(0,0)};
 ?>*\dir{>}\POS?(.5)*+!D{};
 (-25,1.5)*+{e_{-2}};(-15,1.52)*+{e_{-1}};(-5,1.5)*+{e_{0}};(5,1.5)*+{e_{1}};
 (15,1.5)*+{e_{2}};(25,1.52)*+{e_{3}};
 (-25,-4.5)*+{f_{-2}};(-15,-4.5)*+{f_{-1}};(-5,-4.5)*+{f_{0}};
 (5,-4.5)*+{f_{1}};(15,-4.5)*+{f_{2}};(25,-4.5)*+{f_{3}};
 (0,8)*+{g};(0.1,-3)*+{v_0};(10.1,-3)*+{v_1};
 (-9.9,-3)*+{v_{-1}};
 (-19.9,-3)*+{v_{-2}};
 (20.1,-3)*+{v_{2}};
\endxy
\vskip 1pc
\noindent
In fact,
the  elements
$$p_{v_n}:=p_{\{v_n\}},\ s_{e_n}:=p_{\{v_{n-1}\}}s_b,\ s_{f_n}:=p_{\{v_n\}}s_c,\ s_g:=p_{\{v_0\}}s_a,$$
$n\in \mathbb Z$, generates the ideal $I$ and
forms  a Cuntz-Krieger $E$-family.
Since $C^*(E)$ is simple, we have $I\cong C^*(E)$.
\end{remark}

\vskip 1pc

Let $\{s_a, p_A\}$ be a universal representation
of a labelled space $(E,\CL,\CB)$
that generates $C^*(E,\CL,\CB)$.
Then
\begin{eqnarray}\label{nonzeroelements}
 p_{A} s_{\af}\neq 0\ \text{ for } A\in \,\CB\, (A\neq\emptyset),\
 \af \in \CL(AE^{\geq 1}).
\end{eqnarray}
In fact, if $p_{A} s_{\af} =s_\af p_{r(A,\af)}=0$,
then $s_\af^* s_\af  p_{r(A,\af)}=p_{r(A,\af)}=0$,
but $p_{r(A,\af)} $ is nonzero since $r(A,\af)\in \CB$ is non-empty.

\vskip 1pc

\begin{thm}\label{thm-stcofinal} Let $(E, \CL, \bE)$ be a weakly left-resolving
labelled space.
If $C^*(E, \CL, \bE)$ is simple, then $(E, \CL, \bE)$  is strongly cofinal.
\end{thm}
\begin{proof}
Suppose that  $(E, \CL, \bE)$ is not strongly cofinal.
Then   there are $[v]_l$, $x\in  \CL(E^\infty)$,
and  $w\in s(x)$ such that
\begin{eqnarray}\label{assumption-stcofinal}
r([w]_1,  x_{[1,N]})\nsubseteq \cup_{i=1}^m r([v]_l,\ld_i)\end{eqnarray}
for all $N\geq 1$ and any finite number of labelled paths $\ld_1, \dots, \ld_m$.
Consider the ideal
$I$  generated by the projection $p_{[v]_l}$.
 Suppose $ p_{[w]_1}\in I $.
 Then there is an element
$\sum_{j=1}^m c_j(s_{\af_j} p_{A_j} s_{\bt_j}^*)p_{[v]_l}
(s_{\gm_j} p_{B_j} s_{\dt_j}^*)\in I$, $c_j\in \C$, such that
\begin{eqnarray}\label{approx}
\big\|\ \sum_{j=1}^m c_j(s_{\af_j} p_{A_j} s_{\bt_j}^*)p_{[v]_l}(s_{\gm_j} p_{B_j} s_{\dt_j}^*)
 -p_{[w]_1}\ \big\|< 1.
\end{eqnarray}
By Remark~\ref{elements},
we may assume that the paths $\dt_j$'s in (\ref{approx}) have the same length.
Then
\begin{align*}
1>\ & \big\|\ \sum_{j}  c_j(s_{\af_j} p_{A_j} s_{\bt_j}^*)p_{[v]_l}(s_{\gm_j} p_{B_j} s_{\dt_j}^*)
 -p_{[w]_1}\ \big\|\notag\\
\geq\ & \big\|\ \sum_{j} c_j(s_{\af_j} p_{A_j} s_{\bt_j}^*)p_{[v]_l}
(s_{\gm_j}p_{B_j} s_{\dt_j}^*) p_{[w]_1}
 -p_{[w]_1}\ \big\|\notag\\
=\ & \big\|\ \sum_{j} c_j(s_{\af_j} p_{A_j} s_{\bt_j}^*)p_{[v]_l}
(s_{\gm_j} p_{r([v]_l,\gm_j)\cap B_j \cap r([w]_1, \dt_j)}s_{\dt_j}^*)
 -p_{[w]_1}\ \big\|.
\end{align*}

\noindent  We first show that for each $j=1,\dots, m$
\begin{eqnarray}\label{delta-range}
r([w]_1,\dt_j)\subset \cup_{i=1}^m r([v]_l, \gm_i).
\end{eqnarray}
Suppose  $ r([w]_1,\dt_j) \nsubseteq \cup_{i=1}^m r([v]_l, \gm_i)$
for some $j$.
  Then
  $r([w]_1,\dt_j)\setminus \cup_{i=1}^m r([v]_l, \gm_i)\in\bE$ is
  nonempty, hence
  $$ p_j:=p_{r([w]_1,\dt_j)\setminus \cup_{i=1}^m r([v]_l, \gm_i)}\neq 0.$$
  Then with $J:=\{i\mid \dt_i=\dt_j\, \}$,
\begin{align*}
  1> &\ \big\|\ \big(\sum_{i} c_i(s_{\af_i} p_{A_i} s_{\bt_i}^*)p_{[v]_l}
(s_{\gm_i} p_{r([v]_l,\gm_i)\cap B_i \cap r([w]_1, \dt_i)}s_{\dt_i}^*)
 -p_{[w]_1}\big) s_{\dt_j}\ \big\|\\
 = & \ \big\|\  \sum_{i\in J} c_i(s_{\af_i} p_{A_i} s_{\bt_i}^*)p_{[v]_l}
 s_{\gm_i} p_{r([v]_l,\gm_i)\cap B_i \cap r([w]_1, \dt_i)}
 -p_{[w]_1}  s_{\dt_j}\ \big\|\\
  & \ (\text{here we use} \ |\dt_i|=|\dt_j|)\\
 = & \ \big\|\  \sum_{i\in J} c_i(s_{\af_i} p_{A_i} s_{\bt_i}^*)p_{[v]_l}
 s_{\gm_i} p_{r([v]_l,\gm_i)\cap B_i \cap r([w]_1, \dt_i)}
 - s_{\dt_j}p_{r([w]_1,\dt_j)}  \ \big\|\\
  \geq & \ \big\|\  \sum_{i\in J} c_i(s_{\af_i} p_{A_i} s_{\bt_i}^*)p_{[v]_l}
 s_{\gm_i} p_{r([v]_l,\gm_i)\cap B_i \cap r([w]_1, \dt_i)}p_j
 - s_{\dt_j}p_{r([w]_1,\dt_j)} p_j \ \big\|\\
 = & \ \big\|\   s_{\dt_j} p_j \ \big\| =  1,
\end{align*}
which is a contradiction and (\ref{delta-range}) follows.
Also $\dt_j\neq  x_{[1,|\dt_j|]}$ for each $j$.
In fact, if $\dt_j= x_{[1,|\dt_j|]}$,
then by (\ref{delta-range}),
$$r([w]_1,  x_{[1,|\dt_j|]}) =
r([w]_1,  \dt_j) \subseteq \cup_{j=1}^m r([v]_l, \gm_j),$$
which contradicts to (\ref{assumption-stcofinal}).
Thus $s_{\dt_i}^*(s_{x_1} \cdots s_{x_N})=0$ for
$i=1,\dots, m$, and $N=\max_{1\leq j\leq m}\{|\dt_j|\}$.
 Then $y:=p_{[w]_1}s_{x_1} \cdots s_{x_N}$
 is a nonzero partial isometry by (\ref{nonzeroelements}) and
 $s_{\dt_i}^*y=s_{\dt_i}^* p_{[w]_1}s_{x_1} \cdots s_{x_N}
 =p_{r([w]_1,\dt_i)}s_{\dt_i}^*(s_{x_1} \cdots s_{x_N})=0$.
 From (\ref{approx}), we have
\begin{align*}
 1>& \ \big\|\ \big(\sum_{i=1}^m \ld_i(s_{\af_i} p_{A_i} s_{\bt_i}^*)p_{[v]_l}(s_{\gm_i}
p_{B_i} s_{\dt_i}^*)\big)(yy^*)  -p_{[w]_1}(yy^*)\ \big\|\\
 = & \ \|yy^*\|=1,
\end{align*}
a contradiction, and so  $ p_{[w]_1}\notin I $.
\end{proof}

\vskip 1pc

Recall \cite{BP2} that $\af\in \CL(E^{\geq 1})$
with $s(\af)\cap [v]_l\neq\emptyset$
is said to be {\it agreeable}  for $[v]_l$ if $\af=\bt\af'=\af'\gm$ for some
$\af',\, \bt,\, \gm\in \CL(E^{\geq 1})$ with $|\bt|=|\gm|\leq l$.
Otherwise $\af$ is said to be
{\it disagreeable}.
We call $[v]_l$ {\it disagreeable} if there is an $N>0$
such that  for all $n>N$ there is an $\af\in \CL(E^{\geq n})$
that is disagreeable for $[v]_l$.
The labelled space  $(E,\CL, \CE^{0,-})$ is
{\it disagreeable} if for every $v\in E^0$ there is an $L_v>0$
such that $[v]_l$ is disagreeable for all $l>L_v$.
 We say that a labelled space $(E,\CL,\CB)$ is {\it disagreeable} if
 $(E,\CL, \CE^{0,-})$ is disagreeable.
\vskip 1pc

\begin{prop}\label{prop-disagreeable}
Let $(E,\CL, \CB)$ be a labelled space.
Then we have the following.
\begin{enumerate}
\item[(i)] $[v]_l$ is not disagreeable
 if and only if   there is an $N >0$ such that
 every $\af\in \CL(E^{\geq N})$ with $s(\af)\cap [v]_l\neq \emptyset$
is agreeable for $[v]_l$.

\item[(ii)] If $[v]_k$ is not disagreeable,   $[v]_l$ is not disagreeable for all $l>k$.

\item[(iii)]   $(E,\CL, \CB)$ is disagreeable
if and only if
$[v]_l$ is disagreeable for all $l\geq 1$ and  $v\in E^0$.

\end{enumerate}
\end{prop}

\begin{proof}
(i) Note that  $[v]_l$ is  disagreeable if and only if there is a sequence
$(\af_n)_{n=1}^\infty \subset \CL(E^{\geq 1})$,  $|\af_1|<|\af_2|<\cdots$,
of labelled paths that are disagreeable for $[v]_l$.
So $[v]_l$ is not disagreeable if and only if  there is an $N>0$ such that
every labelled path $\af$ with $|\af|\geq N$ and $s(\af)\cap [v]_l\neq\emptyset$
is  agreeable.

(ii)  If  $[v]_k$ is not disagreeable, then there is
$N>0$ satisfying the condition in (i). If
$\af\in \CL(E^{\geq N})$ and $s(\af)\cap [v]_l\neq \emptyset$,
then $s(\af)\cap [v]_k\neq \emptyset$ (since $[v]_l\subset [v]_k$). Hence
$\af=\bt\af'=\af'\gm$ for some $\af',\,\bt,\, \gm\in \CL(E^{\geq 1})$ with $|\bt|=|\gm|\leq k(<l)$,
which means that $[v]_l$ is not disagreeable by (i).

(iii) By definition, $(E,\CL, \CE^{0,-})$ is disagreeable
if $[v]_l$ is disagreeable for all $v\in E^0$ and $l\geq 1$.
For the converse, let $(E,\CL, \CE^{0,-})$ be disagreeable.
Suppose that  $[v]_k$ is not disagreeable for some $v\in E^0$ and $k\geq 1$,
then $[v]_l$ is not disagreeable for all $l>k$ by (ii), a contradiction.
\end{proof}

\vskip 1pc

If $\af$ and $\af'$ are   labelled paths such that
either $\af=\af'$ or $\af=\af'\af''$ for some $\af''\in \CL(E^{\geq 1})$, we call $\af'$ an
{\it initial segment} of $\af$.
If $\bt\in \CL(E^{\geq 1})$,
we write $\bt^\infty$ for the infinite labelled path
$\bt\bt\cdots\in \CL(E^\infty)$.
We call $\bt\in \CL(E^{\geq 1})$   {\it simple}
if there is no labelled path $\dt\in \CL(E^{\geq 1})$ such that
$|\dt|<|\bt|$ and $\bt=\dt^n$ for some $n\geq 1$.

\vskip 1pc
\begin{remark} \label{samelength} If $\gm$ and $\delta$ are simple labelled paths
in $\CL (vE^{\geq 1})$ such that
$\gm^{|\dt|}=\dt^{|\gm|}$, then one can show that
$\gm=\dt$.
\end{remark}

\vskip 1pc

\begin{lem}\label{lemma-disagreeable} Let $(E,\CL, \CB)$ be a
labelled space and $v\in E^0$. If $[v]_l$ is not disagreeable,
there exists an $N> 1$ such that every
 $\af\in \CL([v]_l E^{\geq N})$
 is of the form $\af=\bt^j\bt'$ for some   $\bt\in \CL(E^{\leq l})$
 and an initial segment $\bt'$ of $\bt$.
 Also every $x\in \CL(E^\infty)$ with $s(x)\cap [v]_l\neq \emptyset$ is
 of the form  $x=\bt^\infty$ for a simple labelled path  $\bt\in \CL(E^{\leq l})$.
\end{lem}

\begin{proof}
By Proposition~\ref{prop-disagreeable}(i) there is an $N>1$ such that
every $\af\in \CL(E^{\geq N})$ with $s(\af)\cap [v]_l\neq \emptyset$,
is of the form $\af=\bt\af'=\af'\gm$ for some $\af', \bt, \gm\in \CL(E^{\geq 1})$
with $|\bt|=|\gm|\leq l$.
 If $|\af'|\leq |\bt|$, then  $\bt\af'=\af'\gm$ shows that
 $\af'$ is an initial segment of $\bt$.
 If $|\af'|> |\bt|$, then  $\af=\bt\af'=\af'\gm$ implies that $\af'=\bt\af''$,
 hence  $\af=\bt\af'=\bt^2\af''=\bt\af''\gm$.
 If $|\af''|>|\bt|$, we can repeat the argument until
 we get  $\af=\bt^k \tilde{\af}$ for
 some $\tilde{\af}$ with $|\tilde{\af}|<\bt$.
   But we can always form a (longer) labelled path $\af\sm$ extending $\af$
   (with $|\sm|\geq  |\bt|$) and so
   the above argument shows that $\af\sm=\bt^k\tilde{\af}\sm$ must be of the form $\bt^l\sm'$
   for some $\sm'$  with $|\sm'|<|\bt|$ ($l >k$).
   Hence the path $\tilde{\af}$ is an
   initial segment of $\bt$ and we prove the assertion.

Let $x\in \CL(E^\infty)$ and $s(x)\cap [v]_l\neq \emptyset$.
Then  $x_{[1,N+k]}=\bt_{(k)}^{m_k} \bt_{(k)}'$
for a simple labelled path $\bt_{(k)}\in \CL(E^{\leq l})$
and its initial segment $\bt_{(k)}'$.
Note that if $|\bt_{(k)}|\leq|\bt_{(k')}|$, $\bt_{(k)}$ is an initial path
of  $\bt_{(k')}$.
Since $|\bt_{(k)}|\leq l$ for all $k$, there is a subsequence
$(\bt_{(k_j)})_j$ of $(\bt_{(k)})_k$ such that
$|\bt_{(k_j)}|= |\bt_{(k_1)}|$ for all $j\geq 1$.
Then  with $\bt:=\bt_{(k_1)}$ we finally have $x=\bt^\infty$.
\end{proof}

\vskip 1pc

\noindent
If a labelled space $(E,\CL, \bE)$ is  weakly left-resolving and $v\in E^0$,
then $\{v\}\in \bE$ if and only if $[v]_l=\{v\}$ for some $l\geq 1$.
In Example~\ref{ex-nonsimple}, $\{v\}\in \bE$ for every $v\in E^0$.
Also the following example shows that the condition $\{v\}\in \bE$ for every $v\in E^0$
does not imply that  $\bE=\CEa$.

\vskip 1pc

\begin{ex} In the following labelled graph $(E,\CL)$, we have
 $\{v_i\}\in\bE$ for every $i=1,\dots, 6$, but  $\{v_4\}\notin\CEa$.
\vskip 1pc
\hskip 5pc
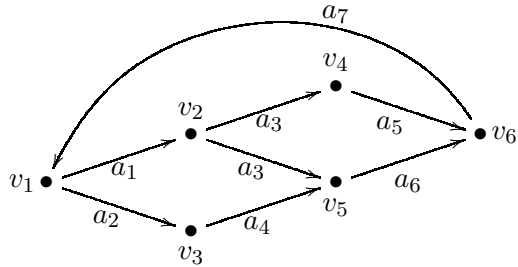
\xy
/r0.38pc/:(0,0)*+{\bullet}="V1";
 (12,4)*+{\bullet}="V2";
 (12,-4)*+{\bullet}="V3";
 "V1";"V2"**\crv{(0,0)&(12,4)};?>*\dir{>}\POS?(.5)*+!D{ };
 "V1";"V3"**\crv{(0,0)&(12,-4)};?>*\dir{>}\POS?(.5)*+!D{ };
 (24,0)*+{\bullet}="V5";
 (24,8)*+{\bullet}="V4";
 "V2";"V5"**\crv{(12,4)&(24,0)};?>*\dir{>}\POS?(.5)*+!D{},
 "V3";"V5"**\crv{(12,-4)&(24,0)};?>*\dir{>}\POS?(.5)*+!D{},
 "V2";"V4"**\crv{(12,4)&(24,8)};?>*\dir{>}\POS?(.5)*+!D{},
 (36,4)*+{\bullet}="V6";
 "V5";"V6"**\crv{(24,0)&(36,4)};?>*\dir{>}\POS?(.5)*+!D{ };
 "V4";"V6"**\crv{(24,8)&(36,4)};?>*\dir{>}\POS?(.5)*+!D{ };
 "V6";"V1"**\crv{(36,4)&(33,10)&(20,15)&(4,10)&(0,0)};?>*\dir{>}\POS?(.5)*+!D{ };
 (6.5,1)*+{a_1},(5,-3)*+{a_2},(18.5,5)*+{a_3},(17,1)*+{a_3},
 (17.5,-3.5)*+{a_4},(28.5,4.8)*+{a_5},(30,0)*+{a_6},(24,14)*+{a_7},
 (-2,0)*+{v_1},(12,6)*+{v_2},(24,10)*+{v_4},(12,-6)*+{v_3},(24,-2)*+{v_5},(38,4)*+{v_6};
\endxy
\end{ex}

\vskip 1pc
For a labelled path $\bt\in \CL(E^{\geq 1})$, put
$T(\bt):=\bt_2\bt_3\cdots\bt_{|\bt|}\bt_1$
whenever $\bt_2\bt_3\cdots\bt_{|\bt|}\bt_1\in \CL(E^{\geq 1})$.
We write $T^k(\bt)$ for $T(T^{k-1}(\bt))$, $k\geq 2$, and
put $T^0(\bt):=\bt$.

\begin{lem} \label{lem-infinite path} Let $(E,\CL, \bE)$ be a
 weakly left-resolving and strongly cofinal labelled space
such that $\{w\}\in \bE$ for every $w\in E^0$.
If $[v]_l$ is not disagreeable,
we have the following.
\begin{enumerate}
\item[(a)]
$\CL(vE^\infty)=\{\bt^\infty\}$ for a simple labelled path $\bt\in \CL(E^{\leq l})$.

\item[(b)]
$|\,r(\{v\},\af)\,|=1 \ \text{ for every }\ \af\in \CL(vE^{\geq 1})$.
\end{enumerate}

\end{lem}
\begin{proof}
(a) By Lemma~\ref{lemma-disagreeable}, every $x\in \CL(vE^\infty)$
is of the form $x=\bt^\infty$
for a simple labelled path $\bt\in \CL(E^{\leq l})$.

\vskip .5pc
\noindent Claim(I)  Let $x,y\in \CL(vE^{\infty})$ and
 $x=\bt^\infty$,  $y=\gm^\infty$
 for simple labelled paths  $\bt,\gm\in \CL(vE^{\leq l})$.
 Then   $\gm=T^k(\bt)$ for some $k\geq 0$
 (thus $|\bt|=|\gm|$).

 \vskip .5pc
 \noindent To prove Claim(I), let $x=x_1x_2\cdots$ and $y=y_1y_2\cdots$,
 $x_i,\,y_i\in \CA$.
 Suppoe there is an $N\geq 1$ such that
 $r(\{v\}, x_{[1,j]})= r(\{v\}, y_{[1,j]})$ for all $j\geq N$.
 Choose $m\geq 1$ with $m|\bt||\gm|>N$ and consider
 the infinite labelled path
 $z:=x_{[\,1,m|\bt||\gm|\,]}y_{[\,m|\bt||\gm|+1,\infty)}\in \CL(vE^\infty)$.
 Then $z$ must be of the form $z=\sm^\infty$
 for a simple labelled path $\sm\in \CL(vE^{\leq l})$ and so
 \begin{eqnarray}\label{z} z =\bt^{m|\gm|} \gm^\infty =\sm^\infty.
 \end{eqnarray}
  Since  $\gm$ and $\sm$ are simple, it follows that
  $\sm =T^k(\gm)$ for some $k\geq 0$.
  Then $|\sm|=|\gm|$ and we have $\bt^{m|\gm|}=\bt^{m|\sm|}=\sm^{m|\bt|}$ from (\ref{z}).
  By Remark~\ref{samelength}, it follows that
  $\bt=\sm=T^k(\gm)$.
  Thus $|\bt|=|\gm|=|\sm|$.
  Then (\ref{z}) shows $\bt=\gm$.
 Thus if $\bt\neq \gm$,  we may assume that there exists a vertex
 $w\in r(\{v\}, y_{[1,j]})\setminus r(\{v\}, x_{[1,j]})$ for some $j$
 large enough with $j> |\bt| |\gm| $.
  Since $(E,\CL, \bE)$ is strongly cofinal,
  there is an $N_1$ and a finite number of labelled paths $\gm_1, \dots, \gm_m$ such that
  $$ r(\{v\}, x_{[1,N_1]})\subset r([v]_1, x_{[1,N_1]})\subset\cup_{i=1}^m r(\{w\}, \gm_i). $$
  If
  $ r(\{v\}, x_{[1,N_1]})\cap r(\{w\},\gm_i)\neq \emptyset$, then
  the labelled path $z:=y_{[1,j]}\gm_i x_{[N_1+1,\infty)}$
  must be of the form
  $z=\sm^\infty$ for a simple labelled path $\sm\in \CL(vE^{\leq l})$.
  Thus
  $$ z=\sm^\infty  =y_{[1,j]}\gm_i  x_{N_1+1}x_{N_1+2}\cdots\bt^\infty,$$
  and so $\sm=T^k(\bt)$ for some $k\geq 0$
  because $\sm$ and $\bt$ are simple.
  Since the initial segment $y_{[1,j]}$ of $z$ has length $j> |\bt| |\gm|$,
  $z$ must be of the form
  $$z=\sm^{|\gm|}\sm^\infty =\gm^{|\sm|}\cdots\sm^\infty,$$
  hence  $\sm^{|\gm|}=\gm^{|\sm|}$.
  Then by Remark~\ref{samelength}, $\sm=\gm$.
  Thus  $\gm= T^k(\bt)$
  for some $k\geq 0$, and Claim(I) is proved.

\vskip 1pc
\noindent Claim(II)  If $x=\bt^\infty,\, y=\gm^\infty\in \CL(vE^{\infty})$
  for simple labelled paths $\bt,\,\gm\in \CL(vE^{\leq l})$,
  then $\bt=\gm$.

 \vskip 1pc
  \noindent Suppose  $\bt\neq \gm$.
  Then by Claim(I), $\gm =T^k(\bt)$ for some $k\geq 1$.
  Let $m=|\bt|=|\gm|$.
  By the first argument in the proof of Claim(I),
  we may assume that there is a vertex
  $u\in r(\{v\}, \gm^n \gm_{[1,j]})\setminus r(\{v\},\bt^n \bt_{[1,j]})$
  for some $n\geq 0$ and $0\leq j\leq n-1$, here $\gm^n \gm_{[1,0]}:=\gm^n $
  and $u$ can be chosen as $u\neq v$.
  By strong cofinality, there exist
  $\dt=\dt_1\cdots\dt_{|\dt|}\in \CL(E^{\geq 1})$   and $N\geq 1$ such that
  $$r(\{v\}, x_{[1,N]})\cap r(\{u\}, \dt)\neq \emptyset$$ and
  \begin{eqnarray}\label{vertex}
  r(\{v\},x_{[1,N-1]})\cap r(\{u\},\dt_{[1,j]})=\emptyset
  \ \ \text{for all } 0\leq j < |\dt|.
  \end{eqnarray}
  (Here $r(\{u\},\dt_{[1,0]}):=\{u\}$.)
  Since $x=\bt^\infty$,
  we can write $$x=x_1 x_2\cdots x_N\cdots=\bt\bt\cdots x_N  \bt_{[i,m]} \bt^\infty$$
  for some $i$.
  Now consider the labelled path
  $$\tilde{y}:=\gm^n \gm_{[1,j]}\dt \bt_{[i,m]}\bt^\infty\in \CL(vE^\infty). $$
  By Claim(I), $\tilde{y}=T^k(\bt)^\infty$ for some $k\geq 1$,
  hence
  $x_N \bt_{[i,m]}\bt^\infty =\dt_{|\dt|} \bt_{[i,m]} \bt^\infty$  and
  we have $x_N=\dt_{|\dt|}$.
  Then
  \begin{align*}
  \emptyset\neq&\  r( \{v\},x_{[1,N]})\cap r(\{u\},\dt)\\
  =&\ r(r(\{v\},x_{[1,N-1]}),x_N)\cap r(r(\{u\},\dt_{[1,|\dt|-1]}),\dt_{|\dt|})\\
  =&\ r(r(\{v\},x_{[1,N-1]}),x_N)\cap  r(r(\{u\},\dt_{[1,|\dt|-1]}),x_N).
  \end{align*}
  But $r(\{v\},x_{[1,N-1]})\cap r(\{u\},\dt_{[1,|\dt|-1]})= \emptyset$
  by (\ref{vertex}),
  a contradiction to that the labelled space is weakly left-resolving,
  and Claim(II) is proved.

\vskip .5pc
\noindent (b) Suppose
$|\,r(\{v\},\af)\,|>1 \ \text{ for an }\ \af\in \CL(vE^{\geq 1})$.
Then there are two paths $\mu,\,\nu\in E^{\geq 1}$ with $s(\mu)=s(\nu)=v$
and $\af =\CL(\mu)=\CL(\nu)$ such that
$v_1:=r(\mu)$ and $v_2:=r(\nu)$ are distinct.
Let $y\in \CL(v_1E^\infty)$.
Since $(E,\CL,\bE)$ is strongly cofinal,
there exist $N\geq 1$ and $\ld_1,\dots,\ld_n\in \CL(E^{\geq 1})$ such that
$$r([v_1]_1,  y_{[1,N]})\subset \cup_{i=1}^n r(\{v_2\}, \ld_i).$$
We can choose  $\ld\in \{\ld_1,\dots,\ld_n\}$ with
$r([v_1]_1,  y_{[1,N]})\cap r(\{v_2\}, \ld)\neq \emptyset$
and may assume that
\begin{eqnarray}\label{single}
r([v_1]_1,  y_{[1,j]})\cap r(\{v_2\}, {\ld}_{[1,i]})= \emptyset
\end{eqnarray}
for all $1\leq j\leq N$ and $0\leq i<|\ld|$.
Since both $\af y$ and $\af\ld y_{[N+1,\infty]}$
belong to $\CL(vE^\infty)$,
by (a)
$$\af y=\af\ld  y_{[N+1,\infty]}=\sm^\infty$$
for a simple labelled path $\sm\in \CL(vE^{\leq l})$.
Thus we have
$$y =\ld  y_{[N+1,\infty]}
= T^k(\sm)^\infty,$$
for some $k$, and we obtain $y_N=\ld_{|\ld|}$.
Note that
\begin{align*}
&\ r(r(\{v_1\}, y_{[1,N-1]}),y_{N})\cap r(r(\{v_2\}, \ld_{|\ld|-1}),\ld_{|\ld|})\\
= & \ r(\{v_1\},y_{[1,N]})\cap r(\{v_2\},\ld)\\
\neq &\ \emptyset\end{align*}
while
$r(\{v_1\}, y_{[1,N-1]})\cap r(\{v_2\}, \ld_{|\ld|-1})=\emptyset$
by (\ref{single}),
a contradiction to the  assumption that $(E,\CL,\bE)$ is weakly left-resolving.
\end{proof}

\vskip 1pc

\begin{thm}\label{thm-disagreeable}
Let $(E,\CL, \bE)$ be a weakly left-resolving labelled space
such that $\{v\}\in \bE$ for each $v\in E^0$.
If $C^*(E,\CL, \bE)$ is simple, then $(E,\CL,\bE)$ is disagreeable.
\end{thm}
\begin{proof}
By Theorem~\ref{thm-stcofinal}, $(E,\CL,\bE)$ is strongly cofinal.
Suppose that $(E,\CL,\bE)$ is not disagreeable.
Then there exists $v\in E^0$ and $l\geq 1$ such that
$[v]_l$ is not disagreeable
by Proposition~\ref{prop-disagreeable}(iii).
Since $\{v\}\in \bE$,
by Proposition~\ref{prop-disagreeable}(ii) and Proposition~\ref{prop-cE}
we may assume that $[v]_l=\{v\}$.
Then, by Lemma~\ref{lem-infinite path},
$\CL(vE^\infty) =\{\bt^\infty\}$ for a simple labelled path
$\bt\in \CL(E^{\leq l})$ and
\begin{eqnarray}\label{singlerange}
|\,r(\{v\},\af)\,|=1 \ \text{ for every }\ \af\in \CL(vE^{\geq 1}).
\end{eqnarray}
Now we consider two possible cases (1) and (2).

\vskip 1pc

\noindent  Case(1)
There is a loop $\mu\in E^{\geq 1}$ at
a vertex $w\in  \{v\}\cup  r(\CL(vE^{\geq 1}))$.
We may assume that $\mu=\mu_1\cdots\mu_{|\mu|}$ is a simple loop,
that is, $r(\mu_i)\neq r(\mu_j)$ for $i\neq j$.
Note from Assumption~\ref{assumptions-graph} and  (\ref{singlerange})
that $\mu$ has no exits and there
are vertices $u_j\in E^0$, $j=1,\dots, |\mu|$,  such that
$$r(\{w\}, \CL(\mu)_{[1,j]})=\{u_j\}.$$
Let $ A:=\{u_1,\dots, u_{|\mu|}\}$.
Then  $A$ and $ \{u_j\}$ belong to $\bE$
so that  the projections
$p_A$ and $p_j:=p_{\{u_j\}}$, $j=1,\dots, |\mu|$, are nonzero
and $p_A$ is the unit of
the $C^*$-subalgebra $p_A C^*(E,\CL,\bE)p_A$ which is simple
as a  hereditary $C^*$-subalgebra of a simple $C^*$-algebra
$C^*(E,\CL,\bE)$.
For  $\gm,\,\dt\in \CL(E^{\geq 1})$, note from (\ref{singlerange}) that
\begin{eqnarray}\label{gm_dt}
r(A,\gm)\cap r(A,\dt)\neq \emptyset
\Longleftrightarrow r(A,\gm)= r(A,\dt)=\{u_j\},
\ j=1,\dots ,|\mu|.\end{eqnarray}
Also for
$s_\gm p_B s_\dt^*\in C^*(E,\CL,\bE)$, $\gm,\dt\in \CL(E^{\geq 1})$
and $B\in \bE$,
if
$p_A(s_\gm p_B s_\dt^*)p_A = s_\gm p_{r(A, \gm)\cap B\cap r(A, \dt)} s_\dt^*\neq 0$,
then
$s_\gm p_{r(A, \gm)} s_\dt^*=s_\gm p_{j} s_\dt^*\neq 0$ for some $j$.
Thus we have
\begin{align*}
      p_A C^*(E,\CL,\bE)p_A
=\ &\,   \overline{span} \{p_A( s_\gm p_B s_\dt^*)p_A
\mid  \gm,\dt\in \CL(E^{\geq 1}),\ B\in \bE\,\}\\
=\ &\,   \overline{span} \{ s_\gm p_j s_\dt^*,
\mid \gm,\dt\in \CL(AE^{\geq 1}u_j),\, j=1,2,\dots,|\mu| \,\}.
\end{align*}
But, since $p_j=\sum_{a\in \CL(u_jE^1) } s_a p_{r(u_j,a)} s_a^*$
 and $ \CL(u_jE^1) =\{\CL(\mu_{j+1})\}$,
\begin{align*}
s_\gm p_j= &\ s_\gm  s_{\CL(\mu_{j+1})}p_{r(u_j,\CL(\mu_{j+1}))}s_{\CL(\mu_{j+1})}^*\\
= &\ s_{\gm \CL(\mu_{j+1})} p_{j+1} s_{\CL(\mu_{j+1})}^*
\end{align*}
for $\gm \in \CL(AE^{\geq 1}u_j)$ and $j=1,\dots, |\mu|$
(here $j+1$ means $1$ if $j=|\mu|$), that is,
$s_\gm p_j \in   p_A C^*(E,\CL,\bE)p_A$.
Also every $\gm\in \CL(AE^{\geq 1}u_j)$ satisfies $\gm=p_{j-1}\gm p_j$,
for some $j=1,\cdots,|\mu|$, where   $p_0:=p_{|\mu|}=p_w$.
Thus $p_A C^*(E,\CL,\bE)p_A$ is the $C^*$-algebra generated by the nonzero partial isometries
$s_j:=p_{j-1}s_{\CL(\mu_j)}p_j$, $j=1,\cdots,|\mu|$ such that
$$s_j^*s_j=s_{j+1}s_{j+1}^*,\ s_i^*s_j=0 \ (i\neq j), \text{ and }
\ \sum_{j=1}^{|\mu|} s_j^*s_j=p_A.$$
Hence  it is a quotient algebra of
$C(\T)\otimes M_{|\mu|}$ which is the graph $C^*$-algebra of the graph with
the vertices $r(\mu_i)$ and the edges $\mu_i$, $i=1,\dots,|\mu|$.
Considering the restriction of the gauge action $\gm_z$, $z\in\mathbb Z$,
on $ C^*(E,\CL,\bE)$ to $p_A C^*(E,\CL,\bE)p_A$,
we see by the gauge invariant uniqueness theorem
(see \cite[Theorem 5.3]{BP1}) that
$p_A C^*(E,\CL,\bE)p_A\cong C(\T)\otimes M_{|\mu|}$,
a contradiction.

\vskip .5pc

\noindent  Case(2)
Suppose that there is no loop at a vertex in
$\{v\}\cup\, r(\CL(vE^{\geq 1})) $.
Recall that
$\CL(vE^\infty) =\{\bt^\infty\}$ for a simple labelled path
$\bt\in \CL(E^{\leq l})$.
 Then $r(\{v\}, \bt^m\bt_{[1,j]})\neq r(\{v\}, \bt^n\bt_{[1,k]})$
for $m\neq n$ or $j\neq k$.
Since $C^*(E,\CL,\bE)$ is simple,
if $I$ denotes the ideal generated by the projection $p_{\{v\}}$,
then $I=C^*(E,\CL,\bE)$.
 Hence there exists $X\in I$ such that
 $\|s_\bt^* s_\bt-X\|<\frac{1}{2}$.
Write
$X =\sum_{i=1}^m \ld_i(s_{\af_i} p_{A_i} s_{\sm_i}^*)p_{\{v\}}(s_{\gm_i} p_{B_i} s_{\dt_i}^*),
\ \ld_i\in \C$.
Since
$$(s_{\af_i} p_{A_i} s_{\sm_i}^*)p_{\{v\}}(s_{\gm_i} p_{B_i} s_{\dt_i}^*)
= s_{\af_i}p_{A_i}p_{r(\{v\}, \sm_i)} s_{\sm_i}^*s_{\gm_i} p_{r(\{v\}, \gm_i)}  p_{B_i} s_{\dt_i}^*,$$
by Lemma~\ref{lemma-disagreeable}
we may assume that $\sm_i$'s and $\gm_i$'s are of the form $\bt^n\bt_{[1,k]}$.
Choose $N_1>0$ large enough so that for every $x\in\CL(E^{\geq N_1})$,
the range vertex set $r(\{v\}, x)$ dose not meet
$r(\{v\}, \sm_i)$ or $r(\{v\}, \gm_i)$ for all $i=1,\dots,m$.
Then with $\{u\}=r(\{v\}, \bt^{N_1})$
(recall that $r(\{v\}, \bt^{N_1})$  is a singleton set),
$$
 p_{\{u\}} s_{\af_i}p_{r(\{v\},\sm_i)}
 = s_{\af_i} p_{r(\{u\},\af_i)}p_{r(\{v\},\sm_i)}
 = s_{\af_i} p_{r(\{v\},\, \bt^{N_1}\af_i)}p_{r(\{v\},\sm_i)}
 =0$$
for all $i=1,\dots,m$ (since $|\bt^{N_1}\af_i|>|\sm_i|$), and so we obtain
$$
\frac{1}{2} >  \ \|p_{\{u\}}(s_\bt^* s_\bt-X)\|
 =  \|p_{\{u\}}p_{r(\bt)}p_{\{u\}}\|
 = \|p_{\{u\}}\|
 =  1,$$
a contradiction.
\end{proof}

\vskip 1pc

\begin{remark}\label{remark-BP2}
Let $\pi_{S,P}$ be a non-zero representation of
$C^*(E,\CL,\CE)$, where $\CE=\bE$ or $\CEa$.
Consider a generalized vertex $[w]_d$ for which
 $P_{[w]_d}$ is a nonzero projection in $C^*(E, \CL,\CE)$.
 Since $[w]_d$ is the disjoint union of a finite number of equivalence classes
 $[w_i]_k$ whenever $k\geq d$, for each $k$ there is an $i$
 such that $P_{[w_i]_k}\neq 0$ as noted in the proof of \cite[Theorem 6.4]{BP2}.
 But it does not mean that we may assume $P_{[w]_d}\neq 0$
 for $d\geq R(w)$ as claimed there.
 For example, consider the labelled graph $C^*$-algebra $C^*(E,\CL,\CE)$ in
 Example~\ref{ex-nonsimple} and
 the ideal $I$ generated by the projection $p_{\{v_0\}}$.
 Let $\pi:C^*(E,\CL,\CE)\to C^*(E,\CL,\CE)/I$ be the quotient map
 (see Remark~\ref{remark-nonsimple}).
 Then
 $$\pi(p_{[v_n]_1})=\pi(p_{E^0\setminus\{v_0\}})\neq 0 \ \text{ for } n\neq 0$$
 but
 $\pi(p_{[v_n]_d})=\pi(p_{\{v_n\}})=0$ for $d\geq R(v_n)(\geq 2)$ and $n\neq 0$.
\end{remark}

\vskip 1pc
\noindent
 Nevertheless,
 if we assume that $(E,\CL,\CE)$ is strongly cofinal and disagreeable,
 a slight modification of the proof of \cite[Theorem 6.4 and Theorem 5.5]{BP2}
 gives the following theorem.

\vskip 1pc

\begin{thm}\label{thm-simple} Let $(E, \CL, \CE)$ be a labelled space that is
strongly cofinal and disagreeable, where $\CE=\bE$ or $\CEa$.
Then $C^*(E, \CL, \CE)$ is simple.
\end{thm}

\vskip 1pc

\begin{cor} Let $(E, \CL, \CE)$ be a labelled space such that
$\{v\}\in \CE$ for each $v\in E^0$, where $\CE=\bE$ or $\CEa$.
Then $C^*(E, \CL, \CE)$ is simple if and only if
$(E, \CL, \CE)$ is  strongly cofinal and disagreeable.
\end{cor}

\end{document}